\documentclass[english,russian]{amsart}
\usepackage{amssymb,amsmath}

\newtheorem{thm}{Theorem}
\newtheorem{crl}[thm]{Corollary}
\newtheorem{prp}[thm]{Proposition}
\newtheorem{lm}[thm]{Lemma}

\newcommand{\der}{\partial}

\newcommand{\sudda}[1]{}

\begin{document}

\title{Logarithmic form of Lagrange inversion formula}

\author{A.S. Dzhumadil'daev}

\address
{Institute of Mathematica, Pushkin street 125,  Almaty,  050000,
Kazakhstan} \email{dzhuma@hotmail.com}


\begin{abstract}
We give presentation of composition inverse of formal power serie in a logarithmic form.
\end{abstract}

\maketitle
\section{Introduction}
For a formal power serie $f(x)=a_1x+a_2\frac{x^2}{2!}+a_3\frac{x^3}{3!}+\cdots\in {\bf C}[[x]]$
its inverse under composition can be given by Lagrange inversion formula
$$f^{\langle -1\rangle}(x)=\sum_{m\ge 1} \left.\frac{d^{m-1}}{dx^{m-1}}\left(\frac{x}{f(x)}\right)\right|_{x=0}\frac{x^m}{m!}$$
(see \cite{Stanley}). The aim of our paper is to prove the following logarithmic version of Lagrange inversion formula

\begin{thm}\label{LogForm}
$$f^{\langle -1\rangle}(x)=\ln \sum_{m\ge 0} \left.\left(\frac{1}{f'(x)}\frac{d}{dx}\right)^m\left(e^x\right)\right|_{x=0}\frac{x^m}{m!}.$$
\end{thm}

To prove this result we introduce three kinds of multiplications on differential operators, composition $\diamond$, white multiplication $\circ$ and black multiplication $\bullet.$ They are partially associative and partially left-symmetric (Proposition \ref{threemult}). 
We give presentation of compositions of differential operators of first order in terms of these multiplications (Theorem \ref{compos}). It allows us to use black multiplication  Bell polynomials to construct powers of differential operators of first order.

\section{Three kinds of multiplications on differential operators}

We consider differential operators of a form $\sum_{i=1}^n u_i\der_i,$ where $u_i=u_i(x_1,\ldots,x_n)$ and $\der_i=\frac{d}{dx_i}$ are partial derivations. Sometimes differential operators of first order are called vector fields.

Let ${\bf Z}_0$ be set of non-negative integers, and 
$${\bf Z}_0^n=\{\alpha=(\alpha_1,\ldots,\alpha_n) | \alpha_i\in {\bf Z}_0\},$$
$$x^{\alpha}=\prod_{i=1}^n x_i^{\alpha_i}, \quad \der^\alpha=\prod_{i=1}^n \der_i.$$
 Denote by $Dif\!f_n$ space of differential operators with $n$ variables. Say that a differential operator $X=\sum_{\alpha\in {\bf Z}}u_\alpha \der^\alpha$ has 
(differential) order  $k$ if $u_\alpha=0$ as soon as $|\alpha|=\sum_{i}\alpha_i\ne k.$ Denote by $Dif\!f_{n,k}$ space of differential operators of order $k.$ In case $k=1$ we use special notation 
$Vect(n)=Dif\!f_{n,1}.$ 

In our paper we study  compositions of differential operators of first order with $n$ variables. We introduce three kinds of multiplications on differential operators.
The first one is a composition, $\diamond$-multiplication,  defined by on generators 
$$u\der^\alpha \diamond v \der^\beta=\sum_{\gamma\in {\bf Z}^n_0} {\alpha\choose \gamma} u\der^\gamma(v)\der^{\alpha+\beta-\gamma}.$$
The second one is $\circ$-multiplication. It is defined by 
$$u\der^\alpha\circ v\der^\beta=u\der^\alpha(v)\der^\beta.$$
The third one is $\bullet$-multiplication. It is defined by 
$$u\der^\alpha \bullet v\der^\beta= u v \der^{\alpha+\beta}.$$
Sometimes we will call $\circ $ and $\bullet$ as white and black multiplications.

\begin{prp}\label{threemult}
The multiplication $\diamond$ is associative. For any $X,Y,Z\in Diff_n,$ 
$$X\diamond(Y\diamond Z)=(X\diamond Y)\diamond Z$$
The multiplication $\bullet $ is associative and commutative. For any $X,Y,Z\in Diff_n,$
$$X\bullet (Y\bullet Z)=(X\bullet Y)\bullet Z),$$
$$X\bullet Y=Y\bullet X.$$
For any $X\in Vect(n)$ and  $Y,Z\in Diff_n,$ 
$$(X,Y,Z)^\circ =(X\bullet Y)\circ Z,$$
$$X\circ (Y\bullet Z)=(X\circ Y)\bullet Z+Y\bullet (X\circ Z),$$
where $ (X,Y,Z)^\circ =X\circ(Y\circ Z)-(X\circ Y)\circ Z$ is associator for the multiplication~$\circ.$ 
For any $X,Y\in Vect(n),$ $Z\in Diff_n,$
$$(X,Y,Z)^\circ =(Y,X,Z)^\circ .$$
\end{prp}

{\bf Proof.} Associativity of the multiplication $\diamond$ is well known. Associativity and commutativity of the black multiplication $\bullet$ is evident. 

Let 
$$X=\sum_{i}u_i\der_i,  Y=\sum_{\alpha} v_{\alpha}\der^\alpha,  Z=\sum_{\beta} w_\beta\der^\beta.$$
Then
$$X\circ (Y\circ Z)=$$
$$\sum_{i,\alpha,\beta} u_i\der_i(v_\alpha\der^\alpha(w_\beta))\der^\beta=$$
$$\sum_{i,\alpha,\beta} u_i\der_i(v_\alpha)\der^\alpha(w_\beta)\der^\beta+
u_iv_\alpha\der^{\epsilon_i+\alpha}(w_\beta)\der^\beta,$$

$$(X\circ Y)\circ Z=$$
$$\sum_{i,\alpha,\beta} u_i\der_i(v_\alpha)\der^\alpha(w_\beta)\der^\beta.$$
Therefore,
$$(X,Y,Z)^\circ=\sum_{i,\alpha,\beta}
u_iv_\alpha\der^{\epsilon_i+\alpha}(w_\beta)\der^\beta=(L\bullet M)\circ R.$$
Another application of this fact. If $X,Y\in Vect(n),$ by commutativity of the multiplication $\bullet,$ 
$$(X,Y,Z)^\circ=(X\bullet Y)\circ Z=(Y\bullet X)\circ Z=(Y,X,Z)^\circ.$$

\begin{crl} \label{rsym1}If $X$ be differential operator of first order, then
$$X\circ(Y\circ Z)=(X\diamond Y)\circ Z.$$
\end{crl}

{\bf Proof.} Note that for differential operator of first order $X$ the composition $X\diamond Y$ can be 
presented as
$$X\diamond Y=X\circ Y+X\bullet Y.$$
Therefore by Proposition \ref{threemult}
$$X\circ (Y\circ Z)=(X\circ Y)\circ Z+(X\bullet Y)\circ Z=(X\diamond Y)\circ Z.$$

\section{Composition of differential operators }

Let $L_k=\sum_{j=1}^n u_{k,j}\der_j,$ $k=1,2,\ldots,m,$  be differential operators of first order. 
Let $[m]=\{1,2,\ldots,m\}$ and $A$ is subset of $[m].$ Suppose that $A=\{i_1,\ldots,i_s\},$ where 
$i_1<i_2<\cdots <i_s.$  Let us denote $i_1=min(A)$ as $h(A)$ and  set $b(A)=\{i_2,\ldots,i_s\}.$  Denote by $L_A^\diamond$ composition of differential operators
$$L_A^\diamond=L_{i_s}\diamond\cdots \diamond L_{i_1}.$$   
Let 
$$L_A=L_{b(A)}^\diamond \circ L_{h(A)}.$$
In other words,
$$L_A=(L_{i_s}\diamond \cdots\diamond  L_{i_2})\circ L_{i_1}.$$

Recall that system of non-empty subsets  ${\pi}=\{A_1,\ldots,A_k\}$ is called  partition of $[m],$ if $ [m]=A_1\cup \cdots \cup A_k,$ and $A_i\cap A_j=\emptyset,$ for $i\ne j.$ Denote by ${\Pi}(m)$ set of partitions of $[m].$ For a partition ${\pi}=\{A_1,\ldots,A_k\}\in {\Pi}(m)$  we set 
$$L_{\pi}^{\bullet}=L_{A_1}\bullet\cdots \bullet  L_{A_k}.$$

\begin{thm} \label{compos}If $L_1,\ldots, L_m$ are differential operators of first order, then 
$$L_m\diamond \cdots\diamond L_1=\sum_{{\pi}\in {\Pi}(m)} L_{\pi}^\bullet.$$
\end{thm}

{\bf Proof.} We use induction on $m.$ For  $m=1$ nothing is to prove. 

Let $A$ be some subset of $[m]$ and $A=\{i_1,\ldots,i_s\},$ such that $i_1<\cdots<i_s.$ 
Let us join element $m+1$ to $A$ and denote obtained set $A'.$ Then  $A'\subseteq [m+1].$ 
Let us prove that 
\begin{equation}\label{equ100}
L_{m+1}\circ L_ A=L_{A'}
\end{equation}
Note that $b(A')=\{i_2,\ldots,i_s,m+1\}=\{m+1\}\cup b(A)$ and $ h(A')=i_1=h(A).$ 
Therefore, by Corollary \ref{rsym1} 
$$L_{m+1}\circ L_ A=L_{m+1}\circ (L_{b(A)}^\diamond  \circ L_{h(A)})=(L_{m+1}\diamond L_{b(A)}^\diamond)\circ L_{h(A)}=L_{b(A')}^\diamond\circ L_{h(A')}=L_{A'}.$$
So, (\ref{equ100}) is established.

Note that partitions of $[m+1]$ can be constructed by partitions of $[m]$ in two ways: either $\{m+1\}$ generates separate block or the element $m+1$ is joined to some block of a partition of $[m].$
For a  partition ${\pi}\in {\Pi}(m+1)$ say that $\pi$ has type $t$ 
if the block that contains $m+1$ has $t$ elements. Denote by ${\Pi}(m+1)^{(1)}$ set of partitions of $[m+1]$ of type $1$ any by ${\Pi}(m+1)^{(>1)}$ set of partitions of type~$t>1.$

Suppose that for $m$ our statement is true,
$$L_m\diamond \cdots\diamond  L_1=\sum_{{\pi}\in {\Pi}(m)} L_{\pi}^\bullet.$$
Let $\pi$ be some partition of $[m]$  and ${\pi}=A_1\cup \cdots\cup A_k.$ 
By (\ref{equ100}) for any $r=1,\ldots,k,$ 
$$L_{m+1}\circ L_{A_r}=L_{A_r'}.$$
Therefore, by Proposition \ref{threemult}
$$L_{m+1}\circ L_{\pi}^\bullet =$$
$$\sum_{r=1}^k L_{A_1}\bullet \cdots L_{A_{r-1}}\bullet L_{A_r'}\bullet
L_{A_{r+1}}\cdots \bullet L_{A_k}=$$
$$\sum_{r=1}^k L_{{\pi'}(r)}^{\bullet},$$
where ${\pi}'(r)$ is a partition of $[m+1]$ constructed by partition ${\pi}=\{A_1,\ldots,A_k\}\in {\Pi}(m),$ by the rule: 
$${\pi}'(r) =B_1\cup\cdots \cup B_k,$$
$$B_1=A_1,\ldots B_{r-1}=A_{r-1}, B_r=\{m+1\}\cup A_r, B_{r+1}=A_{r+1},\ldots, B_k=A_k.$$
Therefore
\begin{equation}\label{equ101}
\sum_{{\pi}\in {\Pi}(m)} L_{m+1}\circ L_{\pi}^\bullet=\sum_{{\mathcal B}\in {\Pi}(m+1)^{(>1)}} L_{\mathcal B}^\bullet
\end{equation}

Note that 
\begin{equation}\label{equ102}
\sum_{{\mathcal B}\in {\Pi}(m+1)^{(1)}} L_{\mathcal B}^\bullet=
\sum_{{\pi}\in {\Pi}(m)} L_{m+1}\bullet L_{\pi}^\bullet
\end{equation}

By Proposition \ref{threemult}
$$L_{m+1}\diamond(L_m\diamond\cdots \diamond L_1)=$$
$$\sum_{{\pi}\in {\Pi}(m)} L_{m+1}\circ L_{\pi}^\bullet+
L_{m+1}\bullet L_{\pi}^\bullet.$$
Hence, by (\ref{equ101}) and (\ref{equ102})
$$L_{m+1}\diamond(L_m\diamond \cdots \diamond  L_1)=\sum_{{\mathcal B}\in {\Pi}(m+1)} L_{\mathcal B}^\bullet.$$
So, inductive step is possible. Our Theorem is proved. 

{\bf Example.} There are $5$ partitions of $[3]=\{1,2,3\},$ 
$${\pi}_1=1-2-3, \quad L_{{\pi}_1}^\bullet=L_3\bullet L_2\bullet L_1,$$
$${\pi}_2=12-3, \quad L_{{\pi}_2}^\bullet=L_3\bullet (L_2\circ L_1),$$
$${\pi}_3=13-2, \quad L_{{\pi}_3}^\bullet=L_2\bullet (L_3\circ L_1),$$
$${\pi}_4=1-23, \quad L_{{\pi}_4}^\bullet=(L_3\circ L_2)\bullet L_1,$$
$${\pi}_5=123, \quad L_{{\pi}_5}^\bullet=(L_3\diamond L_2)\circ L_1.$$
Therefore,
$$L_3\diamond L_2\diamond L_1=L_3\bullet L_2\bullet L_1+L_3\bullet (L_2\circ L_1)+L_2\bullet (L_3\circ L_1)+(L_3\circ L_2)\bullet L_1 +(L_3\diamond L_2)\circ L_1.$$

\section{ Power of vector fields in terms of Bell polynomials}
If otherwise is not stated below we use notation $L^m$ for a  power of differential operator of first order under composition, $L^m=L\diamond \cdots \diamond L.$ 
As far powers of degree $m$ under multiplication $\circ$ and $\bullet$ we use the following notations 
$L^{\circ m}=(\cdots (L\circ L)\cdots)\circ L$ and $L^{\bullet m}=L\bullet  \cdots\bullet L.$

Recall that a sequence of non-decreasing integers  $\lambda=(\lambda_1,\ldots,\lambda_k)$ is called {\it partition} of $m$ and denoted $\lambda\vdash m,$ if $\lambda_1+\cdots+\lambda_k=m.$ 
If $0<\lambda_1\le \lambda_2\le \cdots \le \lambda_k,$ we say that $k=length(\lambda)$ is length of partition $\lambda.$  Suppose that among components of $\lambda$ there are $l_1$ elements $1,$ $l_2$ elements equal to $2,$ etc. There is another way to define partitions. 
A sequence of non-negative integers $(l_1,\ldots, l_m)$ generates partition of $m,$
if $1 \,l_1+2l_2+\cdots+ m\,l_m=m.$
We call $1^{l_1}2^{l_2}\cdots m^{l_m}$ {\it  multiplicity form} or $m$-form of partition $\lambda.$

We recall definition of Bell polynomial  $Y_m(x_1,\ldots,x_m).$ For $\lambda=1^{l_1}\cdots m^{l_m}\vdash m,$ set 
$$k_{l_1,\ldots,l_m}=\frac{m!}{l_1! l_2! \cdots l_m! 1!^{l_1} 2!^{l_2}\cdots m!^{l_m}}.$$
Then 
$$Y_m(x_1,\ldots,x_m)=\sum_{\lambda\vdash m} k_{l_1,\ldots,l_m} x_1^{l_1}\cdots x_m^{l_m}.$$
Bell polynomials are defined over associative commutative algebra with generators $x_1,x_2,\ldots .$ 
In particular we can construct Bell polynomials on differential operators algebra under black multiplication. Instead of $x_i$ we can consider a differential operator of first order $L^{i-1}\circ L.$ 
Let us denote $Y_m^\bullet(L,L\circ L,\ldots, L^{m-1}\circ L)$ Bell polynomial $Y_m(x_1,\ldots,x_m)|_{x_i\rightarrow L^{i-1}\circ L},$ where by associative commutative multiplication we understand the multiplication $\bullet.$ For example,
$$Y_3(x_1,x_2,x_3)=x_1^3+3x_1x_2+x_3,$$
and,
$$Y_3^\bullet(L,L\circ L, L^2\circ L)=L^{\bullet 3}+3 \,L\bullet L^{\circ 2}+L^2\circ L.$$

\begin{thm} \label{L^m one} For any differential operator of first order $L,$ and for any non-negative integer $m,$ 
$$L^m=Y_m^\bullet (L,L\circ L,\ldots,L^{m-1}\circ L).$$
\end{thm}

{\bf Proof.} It is known that number of partitions of $[m]$ as union of $l_i$ subsets with $i$ elements, $i=1,\ldots,m,$ is equal to $k_{l_1,\ldots,l_m}.$ Therefore, by Theorem \ref{compos} 
$$L^m=\sum_{1^{l_1}\cdots m^{l_m}\vdash m}  k_{l_1,\ldots,l_m}\underbrace{\; (L^{1-1}\circ L)\bullet\cdots\bullet (L^{1-1}\circ L)}_{l_1 \mbox{ times }}\bullet\cdots \bullet \underbrace{\; (L^{m-1}\circ L)\bullet \cdots \bullet (L^{m-1}\circ L)}_{l_m \mbox{ times }}.$$
Theorem is proved.

{\bf Example.}
$L^4=L^{\bullet 4}+ 6\, L^{\circ 2}\bullet L^{\bullet 2}+4\,L\bullet (L^2\circ L)+3\, (L^{\circ 2})^{\bullet 2}+L^{3}\circ L.$

Bell polynomials  $Y_m(x_1,\ldots,x_m)$ have many interesting properties. Theorem~\ref{L^m one}
allows us re-write these properties in terms of powers vector fields. Let us give some of such results. 

Exponential generating function for Bell polynomials is 
$$1+\sum_{m\ge 1} Y_m(x_1,\ldots,x_m)\frac{z^m}{m!}=\exp \sum_{m\ge 1} x_m \frac{z^m}{m!}.$$ 
Consider this identity on algebra of differential operators under multplication $\bullet$ and take 
$x_i=L^{i-1}\circ L.$ We obtain the following result.

\begin{thm}\label{exp Lz} For any differential operator of first order $L,$
$$\sum_{m\ge 0} \frac{L^m z^m}{m!}=\exp^\bullet  \sum_{m\ge 1}L^{m-1}\circ L \; \frac{z^m}{m!}$$
\end{thm}

In other words,
$$\exp L\,z=\sum_{k\ge 0} \frac{(\sum_{m\ge 1} L^{m-1}\circ L\,z^m/m!)^{\bullet k}}{k!}.$$
Another formulation of Theorem \ref{exp Lz}
$$ \sum_{m\ge 1}L^{m-1}\circ L \; \frac{z^m}{m!}=\ln^{\bullet}(1+ \sum_{m\ge 1} \frac{L^m z^m}{m!}).$$

{\bf Example.} Let $L=x\der.$ Then $L^{i-1}\circ L=L,$ for any $i\ge 1.$ 
Further, $L^{\bullet k}=x^k\der^k.$
Since Touchard  polynomial $\sum_{k=0}^m S(m,k)x^k$  can be expressed as the value of Bell polynomial on all arguments being $x,$
$$\sum_{k=0}^m S(m,k)x^k=Y_m(x,x,\ldots,x),$$ 
by Theorem \ref{L^m one} we obtain that 
$$(x\der)^m=L^m=Y_m^\bullet (L,L,\ldots,L)=\sum_{k=0}^m S(m,k)L^{\bullet k}=
\sum_{k=0}^m S(m,k)x^k \der^k.$$
where $S(m,k)$ are Stirling numbers of second kind. 
Since
$$\sum_{m\ge 1} L^{m-1}\circ L\frac{z^m}{m!}=\sum_{m\ge 1} L\frac{z^m}{m!}=
L\sum_{m\ge 1}\frac{z^m}{m!}=L(e^z-1),$$
by Theorem \ref{exp Lz} we have also 
$$\exp x\der \, z=\exp^\bullet x\der (e^z-1)=\sum_{i\ge 0} \frac{(x\der)^{\bullet i}(e^z-1)^i}{i!}= \sum_{i\ge 0} \frac{x^i\der^i(e^z-1)^i}{i!}. $$

\section{Logarithmic form for Lagrange inversion formula}

Let ${\bf C}[[x]]$ be algebra of formal power series, ${\bf C}_1[[x]]$ be its subspace generated by formal power series of a form $f(x)=a_1x+a_2\frac{x^2}{2!}+a_3\frac{x^3}{3!}+\cdots,$ where $a_1\ne 0.$ Then ${\bf C}_1[[x]]$ forms a group under composition  $f\diamond g(x)=f(g(x)).$ Unit is the identity function $x.$ Let $f^{\langle -1\rangle} (x)$ be inverse under composition for $f(x)\in {\bf C}_1[[x]],$
$f\diamond f^{\langle -1\rangle} (x)=x.$ Let $f^{\langle -1\rangle} (x)=b_1x+b_2\frac{x^2}{2!}+b_3\frac{x^3}{3!}+\cdots ,$
where $b_1,b_2,b_3,\ldots\in {\bf C}.$
Lagrange inversion formula \cite{Stanley}  states that 
$$b_n=\left.\frac{d^{n-1}}{dx^{n-1}}\left(\frac{x}{f(x)}\right)^n\right|_{x=0},\qquad n\ge 1.$$

Let us construct another form for Lagrange inversion formula. 
\begin{lm} \label{mai1}For any $n\ge 1,$ 
$$b_n=\left.\left(\frac{1}{f'(x)}\frac{d}{dx}\right)^{n-1}\left(\frac{1}{f'(x)}\right)\right|_{x=0}.$$ 
\end{lm}

{\bf Proof.} Let $g(x)=f^{\langle -1\rangle}(x).$ Let us prove that 
\begin{equation}\label{faa}
g^{(n)}(f(x))=\left({\frac{1}{f'(x)}\frac{d}{dx}}\right)^{n-1}\left(\frac{1}{f'(x)}\right).
\end{equation}
By chain rule
$$g(f(x))=x\Rightarrow g'(f(x)f'(x)=1\Rightarrow g'(f(x))=\frac{1}{f'(x)}.$$
So, (\ref{faa}) is true for $n=1.$ Suppose that (\ref{faa}) is valid for $n-1>0.$ By chain rule
$$(g^{(n-1)}(f(x)))'=g^{(n)}(f(x))f'(x).$$
Therefore, 
$$g^{(n-1)}(f(x))=\left({\frac{1}{f'(x)}\frac{d}{dx}}\right)^{n-2}\left(\frac{1}{f'(x)}\right)\Rightarrow$$
$$g^{(n)}(f(x))f'(x)=\left({\frac{1}{f'(x)}\frac{d}{dx}}\right)^{n-2}\left(\frac{1}{f'(x)}\right)'\Rightarrow$$
$$g^{(n)}(f(x))=\left({\frac{1}{f'(x)}\frac{d}{dx}}\right)^{n-1}\left(\frac{1}{f'(x)}\right)$$
So, induction by $n$ is valid, and  (\ref{faa}) is established. 

Since $f(0)=0,$ by (\ref{faa})
$$g^{(n)}(0)=\left.\left({\frac{1}{f'(x)}\frac{d}{dx}}\right)^{n-1}\left(\frac{1}{f'(x)}\right)\right|_{x=0}.$$
Lemma  \ref{mai1} is proved.

{\bf Proof of Theorem \ref{LogForm}.} Follows from Theorem \ref{exp Lz} and Lemma \ref{mai1}.

{\bf Example.} Let $f(x)=xe^{-x}.$ Then (\cite{Stanley}, Example 5.4.4)
$$f^{\langle -1\rangle}(x)=\sum_{m\ge 1}m^{m-1}\frac{x^m}{m!}.$$
Note that $f'(x)=e^{-x}(1-x)$ and
$$\left.\left(\frac{e^x}{1-x}\,\frac{d}{dx}\right)^m\left(e^x\right)\right|_{x=0}=(m+1)^{m-1}.$$
Therefore, 
$$f^{\langle -1\rangle}(x)=\ln \sum_{m\ge 0} (m+1)^{m-1}\frac{x^m}{m!}.$$

\end{document}